\documentclass[leqno, 12pt]{amsart}
\usepackage{amssymb,amsmath,amsthm,graphicx,multirow,setspace,url}
\usepackage{amscd}
\usepackage{tikz-cd}

\usepackage{graphicx}

\usetikzlibrary{arrows}
\usetikzlibrary{quotes}

\title{Local Cohomology of Module of Differentials of integral extensions}
\author{S. P. Dutta}
\address{Department of Mathematics\\
University of Illinois at Urbana-Champaign\\
1409 West Green Street\\
Urbana, Illinois 61801}

\begin{document}

\newcommand{\hookuparrow}{\mathrel{\rotatebox[origin=c]{90}{$\hookrightarrow$}}}
\newcommand{\hookdownarrow}{\mathrel{\rotatebox[origin=c]{-90}{$\hookrightarrow$}}}


\theoremstyle{plain}
\newtheorem*{theorem*}{Theorem}
\newtheorem{theorem}{Theorem}

\newtheorem{innercustomthm}{Theorem}
\newenvironment{customthm}[1]
  {\renewcommand\theinnercustomthm{#1}\innercustomthm}
  {\endinnercustomthm}

%



\newtheorem*{corollary}{Corollary}

\theoremstyle{remark}
\newtheorem*{remark}{Remark}
\newtheorem*{claim}{Claim}

\newtheorem*{namedtheorem}{\theoremname}
\newcommand{\theoremname}{testing}
\newenvironment{named}[1]{\renewcommand{\theoremname}{#1}
  \begin{namedtheorem}}
	{\end{namedtheorem}}

\newcommand{\gr}{\operatorname{grade}}
\newcommand{\Syz}{\operatorname{Syz}}
\newcommand{\syz}[2]{\Syz^{#1}(#2)}
\newcommand{\ring}[2]{{\mathcal{O}_{#1}(#2)}}
\newcommand{\hm}[3]{H_{#1}^{#2}(#3)}
\newcommand{\Tor}{\operatorname{Tor}}
\newcommand{\tor}[3]{\Tor_{#1}^{#2}(#3)}
\newcommand{\Ext}{\operatorname{Ext}}
\newcommand{\ext}[3]{\Ext_{#1}^{#2}(#3)}
\newcommand{\Id}{\operatorname{Id}}
\newcommand{\im}{\operatorname{Im}}
\newcommand{\coker}{\operatorname{coker}}
\newcommand{\grade}{\operatorname{grade}}
\newcommand{\Ht}{\operatorname{height}}
\newcommand{\Hom}{\operatorname{Hom}}
\newcommand{\Der}{\operatorname{Der}}
\newcommand{\dm}{\operatorname{dim}}
\renewcommand{\H}{\text{H}}
\newcommand{\Tr}{\text{Tr}}
\newcommand{\D}{\text{D}}
\renewcommand{\L}{\text{L}}

\newcommand{\ul}[1]{\underline{#1}}

\begin{abstract}
The main focus of this paper is on determining the highest non-vanishing local cohomology modules of $\Omega_{B/R}, \Omega_{B/V} \\ (\Omega_{B/k})$ where $R$ is either a complete regular local ring or a complete local normal domain with coefficient ring $V$(field $k$) and $B$ is its integral closure in an algebraic extension of $Q(R)$. Similar problem is also studied over a normal domain $R$ containing a field k of characteristic 0. In this connection new observations on the direct summand property for integral extensions are also presented.
\end{abstract}

\maketitle

\section{Introduction}\label{Intro}
In this note we mainly address the following problem: let $R$ be a normal domain containing a field k of characteristic 0 or let $(R, m)$ be a complete regular local ring/complete local normal domain of dimension $n$ with co-efficient ring $V$ (field $k$). Let $K$ denote the field of fractions of $R$ and let $B$ denote the integral closure of $R$ in an algebraic extension $F$ of $K$. Let $\Omega_{B/R}, \Omega_{B/V} (\Omega_{B/k})$ denote the module of differentials of $B$ over $R, V(k)$ respectively.
\smallskip

\textbf{Problem.} Determine whether $\Hom_R (\Omega_{B/k}, R) (\Hom_R (\Omega_{B/V}, R)) \neq 0$, $\Hom_R (\Omega_{B/R}, R)) \neq 0$ or determine the highest non-vanishing local cohomology modules with respect to $m$ of $\Omega_{B/V} (\Omega_{B/k})$ and $\Omega_{B/R}$.

In section 2 we prove several lemmas that will be used in the proofs of theorems in section 3 \& section 4. The corresponding results include providing sufficient conditions for a) extensions $A \rightarrow B$ so that $\Omega_{B/A}$ becomes $B$-flat, $\Gamma_{B/A}$ becomes null and for b) extensions $A \rightarrow B \rightarrow C$ so that the corresponding right exact sequence of module of differentials becomes left exact. Another useful tool in this respect is theorem 2.3 where by utilizing Popescu's work ([11], [14]) we provide a simple proof of the fact that for any geometrically regular extension of Noetherian rings $A \rightarrow B, \Omega_{B/A}$ is $B$-flat and $\Gamma_{B/A} = 0$.

\line(1,0){300}

AMS Subject Classification:
13D02, 13D22, 13D45, 13N05

Key words and phrases: module of differentials, integral extension, geometric regularity, local cohomology, formally unramified extensions.

 The second part of this result is originally due to M. Andr\'e ([1]) and later it was proved by Popescu ([11]) via his characterization of
geometrically regular extensions. A stronger version of this result has been proved by Gabber \& Ramero (Cor. 5.6.25, [4]).

As a corollary to lemma 2.1 it has been pointed out that for any integral domain $A$ integral over a complete local normal domain $(R, m)$ of dimension $n$ such that $Q(A)$ is separably algebraic over $Q(R)$, $\H_m^n (\Omega_{A/R}) \\ =0$.

One of our main results in section 3 is the following.

\smallskip
\textbf{Theorem 3.5.} Let $R$ be a normal domain of dimension $n$ containing a field $k$ of characteristic 0 with a non-null derivation $D \ \in \Der_k (R)$ and let $B$ be its integral closure in an algebraic extension $F$ of $K$. Then $\Hom_R(\Omega_{B/k}, R) \neq 0$. If $(R, m)$ is a complete local normal domain of dimension $n$ with maximal ideal $m$ and a non-null derivation $D \in \Der_k (R)$, then $\H{_m^n}(\Omega_{B/k})\neq 0$.

The two key steps in the proof of this theorem are \textbf{Theorem 3.1} and \textbf{Proposition 3.4} where the same statement is proved for a complete d.v.r. and for a d.v.r. with a non-null derivation respectively. In these proofs we provide a prescription for defining a non-null $R$-linear map: $\Omega_{B/k} \rightarrow R$ that depends on the non-null derivation $D$.

\smallskip
Next we deal with complete local domains. We would like to mention here that if $(C, q)$ is a complete local domain  of dimension $n$ with coefficient ring $V(k)$, then $\Omega_{C/V} (\Omega_{C/k})$ is not necessarily a finitely generated $C$-module. We prove the following:

\smallskip
\textbf{Theorem 3.6.} Let $(C, q)$ be a complete local domain of dimension $n \geq 2$ with coefficient ring $V$ (field $k$). Let $(R, m)$ be a power series ring over $V(k)$ contained in $C$ such that $C$ is a module finite extension of $R$. We have the following: \\
i) if $Q(C)/Q(R)$ is separably algebraic then $\H_q^n(\Omega_{C/V}) (\H_q^n(\Omega_{C/k})) \neq 0$, \\
ii) if characteristic of $k =p >0, k$ perfect and $\Omega_{Q(C)/k} \neq 0$, then $\H_q^n(\Omega_{C/k}) \neq 0$ and  \\
iii) if $Q(C)/Q(R)$ is separably algebraic, then $\H_q^n(\Omega_{C/R}) = 0$; moreover if $\Omega_{C/R} \neq 0$ and $C$ is normal, then $\H_q^{n-1}(\Omega_{C/R}) \neq 0$.

\smallskip
\textbf{Theorem 3.7.}, the final result of this section, can be viewed as an extension of the theorem on purity of branch locus ([10], [15]) to infinitely generated integrally closed domains $B$ integral over a regular local ring $R$.

\smallskip
In section 4 we deal with the situation where $R$ is a complete regular local ring and $B$ is the absolute integral closure of $R$. Our problem is irrelevant in positive characteristic, since $\Omega_{B/k} = 0$ and $\Omega_{B/R} = 0$ i.e. $k \rightarrow B$ and $R \rightarrow B$ are unramified extensions. Due to theorem 3.5 we focus our attention on mixed characteristic. In this part of our work the highest non-vanishing local cohomology of $\Omega_{B/R}$ gets connected with the direct summand property for integral extensions of $R$ conjectured by M. Hochster ([7]). Hochster proved the equicharacteristic case of his conjecture in early seventies ([7]). In 2016 Y. Andr\'e ([2]) proved this conjecture in mixed characteristic by using aspects of Almost Ring Theory introduced by Faltings in his work on $p$-adic Hodge theory. Afterwards Bhatt ([3]) gave a short proof of this conjecture and its derived variant. Both proofs use the notion of perfectoid geometry. The main theorem in section 4 is the following:

\smallskip

\textbf{Theorem 4.1.} Let $(R, m)$ be a complete regular local ring of dimension n with coefficient ring $V$ of mixed characteristic $p$ and let $B$ denote its absolute integral closure. We have the following:

{\bf i)} $\Omega_{B/V} = p \ \Omega_{B/V}, \Omega_{B/R} = p \ \Omega_{B/R}$. These observations imply the following:

\smallskip
a) For any ideal $I$ of $R$ of the form $I=(p, x_1, .., x_{i-1})$ generated by a part of a system of parameters containing $p$ of $R$, $H_I^i(\Omega_{B/V})=0$ and $H_I^i(\Omega_{B/R})=0$. In particular $H_m^n(\Omega_{B/V})=0$ and $H_m^n(\Omega_{B/R})=0$. 

\smallskip
b) $\widehat{\Omega}_{B/V}=0=\widehat{\Omega}_{B/R}$, where  $(\widehat{-})$  denotes the I-adic completion. And hence $V \rightarrow B, R \rightarrow B, V \rightarrow \widehat{B}$ and $R \rightarrow \widehat{B}$ are all formally unramified extensions.

{\bf ii)} Non-vanishing of $\H{_m^{n-1}}(\Omega_{B/R})$ implies the direct summand property for integral extensions of $R$ (equivalently descent of flatness for integral extensions of Noetherian rings ([12])) and the converse is also true.

Non-vanishing of $\H{_m^{n-1}} (\Omega_{B/V})$ implies the direct summand property for integral extensions over $R$.

\smallskip
For part i) we do not need $R$ to be regular. Actually we prove a much more general statement.

\smallskip
As a \textbf{Corollary}, due to Y. Andr\'e's proof of the direct summand conjecture in mixed characteristic, it follows that $\H{_m^{n-1}} (\Omega_{B/R}) \neq 0$.

\smallskip
\ \ \ \ We do not know whether part ii) of the above theorem is valid in equicharacteristic zero.

\smallskip
As corollaries to \textbf{Proposition 4.2} we derive the following over $R=V [[X_1, \dots, X_{n-1}]]$ and $m, B$ as above: i) for any $i>0$, if $H_m^{i+1} (B) \neq 0$, then $H_m^i (\Omega_{B/R}) \neq 0$ and ii) for any $i \geq 0$ there exists an exact sequence
$$0 \rightarrow \underset{1}{\overset{n-1}{\oplus}} \Ext {^{n-i-1}_R} (B, R) \rightarrow \Ext{_R^{n-i}} (\Omega_{B/R}, R) \rightarrow \Ext{_R^{n-i}} (\Omega_{B/V}, R) \rightarrow 0$$

In particular
$$0 \rightarrow \underset{1}{\overset{n-1}{\oplus}} B^* \rightarrow \Ext{_R^1} (\Omega_{B/R}, R) \rightarrow \Ext{_R^1} (\Omega_{B/V}, R) \rightarrow 0$$
is exact.

\smallskip
\ \ \ \ \ \ \ \ Our final theorem deals with local cohomology with respect to the ideal $(x_1, \dots, x_{n-1})$ generated by a part of a regular system of parameters. The second assertion in this theorem provides an apparently approachable sufficient condition for the direct summand property.
\smallskip

\textbf{Theorem 4.3.} Let $(R, m)$ be a complete regular local ring in mixed characteristics $p > 0$ with co-efficient ring $V$, and $B$ be its absolute integral closure. Let $p, x_1, \dots, x_{n-1}$ be a system of parameters of $R$ such that $x_1, \dots, x_{n-1}$ form a regular system of parameters. Let $\underline{x}$ denote the ideal $(x_1, \dots, x_{n-1})$ in $R$. Then $\H{_{\underline{x}}^{n-1}}(\Omega_{B/V}) \neq 0$. Non-vanishing of $\H{_{\underline{x}}^{n-1}}(\Omega_{B/R})$ implies the direct summand property for integral extensions of $R$.

\smallskip

\textbf{Corollary.} If $\H{_p^0} (\H{_{\underline{x}}^{n-1}} (\Omega_{B/V})) \neq 0$, then $\H{_m^{n-1}}(\Omega_{B/V}) \neq 0$.

\smallskip
Let us note that some of the results proved in section 4 are not valid in equicharacteristic zero. I do not know answers to several questions that came up during the course of this work. A few of these questions related to this paper are mentioned at the end.

\smallskip
\ \ \ \ \ \ For definition of $\Gamma_{S/R}$ and for basic properties of $\Omega_{S/R}$ and $\Gamma_{S/R}$ the reader is referred to ([5], [14]).

\smallskip
\textbf{Notations.} Given an extension $R \rightarrow S$ the corresponding module of differentials is denoted by $\Omega_{S/R}$ (instead of $\Omega^1_{S/R}$); d.v.r. stands for discrete valuation ring; for any integral domain (ring) $C, Q(C)$ denotes the field of fractions (quotient ring) of $C$; for any local ring $(R, m), E$ denotes the injective hull of $R/m$ over $R$, $E(k (P))$ denotes the injective hull of $R/P$ over $R$ and for any $R$-module $M, M^{\vee}$ stands $\Hom_R (M, E)$, $M^*$ stands for $\Hom_R (M, R)$.

\smallskip
\section*{Section 2.} \label{s2} In this section we are going to prove several lemmas and propositions that will be used in section 3 and section 4.

\smallskip
\textbf{2.1. Lemma.} Let $R$ be a normal domain, $R \hookrightarrow A$ an integral extension such that $A$ is an integral domain and  $Q(R)\hookrightarrow Q(A)$ is a separable algebraic extension. Then $\Hom_R(\Omega_{A/R}, R)=0$.

\smallskip
\textbf{Proof.} For any $a \in A$, $a \neq 0$, the minimal polynomial $f(X)$ for $a$ has all its co-efficient in $R$; separability of $Q(A)$ over $Q(R)$ implies $f'(a) \neq 0$. Since $f(a) = 0, f'(a)da = 0$ \dotfill(1)

Since $f'(a)$ is integral over $R$ and $A$ is an integral domain, an integral equation of smallest degree satisfied by $f'(a)$ is of the form:
$$[f'(a)]^t+r_1[f'(a)]^{t-1}+\dots+r_t=0, r_i \in R$$
where $r_t \neq 0$. This implies, due to (1), that $r_t da = 0$. This is true for every $a$ in $A-R$. Thus $\Omega_{A/R}$ is a torsion module over both $A$ and $R$ and hence $\Hom_R(\Omega_{A/R}, R) = 0$

\medskip
\textbf{Corollary.} If $(R, m)$ is a complete local normal domain of dimension $n$ and $A$ is as above then $\H{_m^{n}}(\Omega_{A/R}) = 0$.

This follows from local duality due to the fact that the canonical module of $R$ is isomorphic to a height 1 ideal in $R$ and $R$ can be embedded into this ideal.

\smallskip
\textbf{2.2. Sub-Lemma.} Let $R \hookrightarrow S$ be two integral domains such that $S = R[X_1, \dots, X_n]/(f_1, \dots, f_n)$, where $f_1, \dots, f_n$ form an $R[X_1, \dots, X_n]$-sequence and $Q(S)/Q(R)$ is separably algebraic. Then $\Gamma_{S/R} = 0$.

\smallskip
\textbf{Proof.} Let $I = (f_1, \dots, f_n)$. We have the following exact sequence of $S$-modules:

$0 \rightarrow \Gamma_{S/R} \rightarrow I/I^2 \rightarrow \oplus{_1^n} Sdx_i \rightarrow \Omega_{S/R} \rightarrow 0$ \dotfill (2)

Since $I/I^2$ is a free $S$ module of rank $n$ and $S$ is a domain, $I/I^2$ is torsion-free over $S$; hence so is $\Gamma_{S/R}$ unless it is 0.

Since $Q(S)/Q(R)$ is separably algebraic, $\Omega_{Q(S)/Q(R)} = 0$, i.e. $\Omega_{S/R} \otimes Q(S) = 0$. Tensoring the above sequence with $Q(S)$ we obtain the following exact sequence:

$0\rightarrow Q(S) \otimes \Gamma_{S/R} \rightarrow Q(S) \otimes I/I^2 \rightarrow \oplus{_1^n} Q(S)dx_i \rightarrow 0$  \dotfill  (3)

Since rank $Q(S) \otimes I/I^2 = n =$ rank $\oplus{_1^n} Q(S)dx_i$, it follows from (3) that $Q(S) \otimes \Gamma_{S/R} = 0$. Hence $\Gamma_{S/R} = 0$.

\smallskip
\textbf{Lemma.} Let $R \hookrightarrow S$ be two regular rings such that $S$ is a module-finite extension of $R$ and $Q(S)$ is separable algebraic over $Q(R)$. Then $\Gamma_{S/R} = 0$.

\smallskip
\textbf{Proof.} First we assume that both $R$ and $S$ are local.

Since dim$R =$ dim$S$, we have $S \simeq R[X_1, \dots, X_n]/P$ where $P$ is a prime ideal of height $n$ (one could also write $S = R[X_1, \dots, X_n]_T/P$ where $T$ is the inverse image of the maximal ideal of $S$ under the onto map $R[X_1, \dots, X_n] \rightarrow S$). Since $S$ is regular local, $P$ is a complete intersection ideal generated by $n$ elements and $P/P^2$ is a free $S$-module of rank $n$. The proof now follows by the above sub-lemma.

Now let us assume that $R$ is local. Let $m$ denote the maximal ideal of $R$ and let $m_1, \dots, m_r$ denote the maximal ideals of $S$. Let $S_i = S_{m_i}$ for $1 \leq i \leq n$ and let $\widehat{R}, \widehat{S}$ and $\widehat{S_i}$ denote the $m$-adic completions of $R,S$ and $S_i$ respectively. Then $\widehat{S_i}$ is a complete regular local ring module-finite over $\widehat{R}$. It can be checked that $\Gamma_{S/R} \otimes \widehat{S} \simeq  \Gamma_{\widehat{S}/\widehat{R}} \simeq {\underset{i=1}{\oplus^r}}  \ \Gamma_{\widehat{S_i}/\widehat{R}}$. Due to the result proved above it follows that $\Gamma_{S/R} \otimes \widehat{S} = 0$ and hence $\Gamma_{S/R} = 0$. \

\ \ \ \ \ \ Finally the general case follows via localization.

\smallskip
\textbf{2.3. Theorem ( Andr\'e, Gabber-Ramero, Popescu,} see introduction \textbf{)} Let $A \rightarrow B$ be an extension of Noetherian rings such that $B$ is geometrically regular over $A$. Then $\Omega_{B/A}$ is $B$-flat and $\Gamma_{B/A} = 0$.

\smallskip
\textbf{Proof.} Popescu's theorem ([11], [14]) asserts that in the above situation $B = \underrightarrow{\lim} \ C$, a filtered inductive limit over smooth $A$-algebras $C$. Then $\Omega_{C/A}$ is a finitely generated projective C-module and $\Gamma_{C/A} = 0$. Since $\Omega_{B/A} = \underrightarrow{\lim} \ \Omega_{C/A}$ $ \Gamma_{B/A} = \underrightarrow{\lim} \ \Gamma_{C/A}$, the conclusion follows.

\smallskip
\textbf{2.4. Lemma.}
$A \rightarrow B \rightarrow C$ are injective ring homomorphisms of integral domains such that $Q(C)$ is a separable algebraic extension of $Q(B)$ and $\Omega_{B/A}$ is a flat $B$-module. Then the following sequence:
$$ 0 \rightarrow \Omega_{B/A} \otimes_BC \rightarrow \Omega_{C/A} \rightarrow \Omega_{C/B} \rightarrow 0$$
is exact.
\smallskip

\textbf{Proof.} we have an exact sequence:

$0 \rightarrow \Gamma_{C/B/A} \rightarrow \Omega_{B/A} \otimes_BC \rightarrow \Omega_{C/A} \rightarrow \Omega_{C/B} \rightarrow 0$ \dotfill (4)

Since $\Omega_{B/A}$ is $B$-flat, $\Omega_{B/A} \otimes_BC$ is $C$-flat; since $C$ is an integral domain this implies that $\Gamma_{C/B/A}$ is a torsion-free $C$-module unless it is null.
Let $S=C-\{0\}$. Tensoring (4) with $S^{-1}C=Q(C)$ we obtain the following exact sequence:
$$0 \rightarrow S^{-1} \Gamma_{C/B/A} \rightarrow \Omega_{Q(B)/Q(A)} \otimes Q(C) \rightarrow \Omega_{Q(C)/Q(A)}\rightarrow \Omega_{Q(C)/Q(B)} \rightarrow 0$$
Since $Q(C)$ is separable algebraic over $Q(B)$, $Q(C)$ is quasi-smooth (in Grothendieck's language: formally smooth with respect to discreet topology) over $Q(B)$ and hence it follows that $S^{-1} \Gamma_{C/B/A} = 0$ (20.5.7, [5]). Since $\Gamma_{C/B/A}$ is torsion-free, it follows that $\Gamma_{C/B/A} = 0$.

\smallskip
\textbf{2.5.} Here we state a result that has been used several times in this paper. This result is known as Jacobi-Zariski sequence, henceforth JZ sequence.

\smallskip
\textbf{JZ sequence.}([1], [14]) Given ring extensions $A \rightarrow B \rightarrow C$, there exists an exact sequence $\Gamma_{C/A} \rightarrow \Gamma_{C/B} \rightarrow \Omega_{B/A} \otimes C \rightarrow \Omega_{C/A} \rightarrow \Omega_{C/B} \rightarrow 0.$ If $\Omega_{B/A}$ is $B$-flat then this sequence can be extended by putting $\Gamma_{B/A} \otimes C$ on the left.

\ \ \ \ \ \ For a proof we refer the reader to ([14]).

\section*{Section 3: Equicharacteristic 0, complete local domains and extension of purity} \label{s2}

\textbf{3.1. Theorem.}
Let $(R, m)$ be a complete discrete valuation ring containing a field $k$ of characteristic 0 and let $B$ denote its integral closure in an algebraic extension $F$ of $Q(R)$. Then $\Hom_R(\Omega_{B/k}, R) \neq 0$.

\textbf{Proof.} Let $K$ denote the field of fractions of $R$. Since characteristic of $k=0$, $R$ is geometrically regular over $k$. By theorem 2.3, $\Omega_{R/k}$ is a flat $R$-module and hence torsion-free over $R$. Let $R \subset C \subset C' \subset B$ be such that $C, C'$ are integral closures of $R$ in $Q(C), Q(C')$ respectively and $[Q(C): K] < \infty$, $[Q(C') : K] < \infty$. Then $C, C'$ are module-finite extensions of $R$ and $C'$ is a module finite extension of $C$. Since $R$ is a complete d.v.r., $C, C'$ are also the same; hence, by argument as above, $\Omega_{C/k}(\Omega_{C'/k})$ is a flat $C(C')$ module and $\Gamma_{C'/C}=0$ (lemma 2.2). Due to lemma 2.4 and the JZ sequence 2.5 we obtain the following short exact sequences:

\smallskip
$0 \rightarrow \Omega_{C/k} \otimes B \rightarrow \Omega_{B/k} \rightarrow \Omega_{B/C} \rightarrow 0$   \dotfill (5)

\smallskip
$0 \rightarrow \Omega_{R/k} \otimes C \rightarrow \Omega_{C/k} \rightarrow \Omega_{C/R} \rightarrow 0$   \dotfill (6)

\smallskip
$0 \rightarrow \Omega_{C/k} \otimes C' \rightarrow \Omega_{C'/k} \rightarrow \Omega_{C'/C} \rightarrow 0$   \dotfill (7)

\smallskip
$0 \rightarrow \Omega_{C/R} \otimes C' \rightarrow \Omega_{C'/R} \rightarrow \Omega_{C'/C} \rightarrow 0$   \dotfill (8)

\smallskip

Also note that
$\Omega_{B/k} = \underrightarrow{\lim} \ \Omega_{C/k}$ and $\Omega_{B/R} = \underrightarrow{\lim}  \ \Omega_{C/R}$ over all C as above \ \ \ \         \dotfill (9)

\medskip
Since characteristic of $k = 0$ and $R$ is a complete d.v.r., any $C$ as above is of the form $C = R [X]/(f(X)) = R[b]$ where $f(X)$ is monic irreducible in $R[X]$ ([13]). Thus $C$ is a free $R$-module with bases $1, b, b^2, \dots, b^{n-1}, n=\deg f(X)$. Hence $\Omega_{C/R} \simeq C/(f'(b)) db$ is a cyclic $C$-module generated by $db$. It follows from (5), (6), (7) and (9) that any $\omega  \in  \Omega_{B/k}$ can be expressed as:

\smallskip
$\omega = \Sigma t_i \ \omega_i + cdb, \ c \in C$ for some $C$ as above, $\omega_i  \in  \Omega_{R/k} \ t_i \in C, \\ 1 \leq i \leq r.$   \dotfill (10)

\smallskip
Let $\Tr: F \rightarrow{K}$ denote the trace-map: for any $x \neq 0 \in F$, if $x \in \L \subset F$ such that $[\L:K] < \infty$, then $\Tr(x) = \Tr_{\L/K}(x)/[\L:K]$. Let $\D  \in  \Der_k(R)$, be non-null (e.g. $D=\partial/\partial X, R = \Bbbk[[X]], \Bbbk = R/m)$.
Since characteristic of $k=0$, $\D$ can be extended uniquely to a derivation $\D$ (same notation): $B \rightarrow F$. We define an $R$-linear map $\L:\Omega_{B/k} \rightarrow K$ as follows: for any $\omega \ \in  \Omega_{B/k}, \L(\omega) = \Tr(\tilde{\D}(\omega))$, where $\tilde{\D}:\Omega_{B/k} \rightarrow F$ is the $B$-linear map induced by $\D$.

\smallskip
{\bf Claim.} Im $\L \subset R$

\smallskip
{\bf Proof of the Claim.} Due to (10), given any $\omega  \in  \Omega_{B/k}$,  we have $\L (\omega) = \L (\Sigma t_i \omega_i) + \L(cdb)$, For any $x  \in  B,  \Tr (x)  \in  R$ and $\tilde{\D}(\omega_i)  \in  R$ for every $i$.

Hence $\L(\Sigma t_i \omega_i)  \in  R$       \dotfill (11).

We need to show that $\L(cdb)  \in  R$. Recall that $f(b) = 0$ where $f(X)$ is an irreducible polynomial in $R[X]$. Let $f(X) = X^n + r_1 X^{n-1} + \dots + r_n, r_i  \in  R$. Then we have $f'(b)db + \Sigma{^n_1} d(r_i)b^{n-i} = 0$.

Hence $f'(b) \tilde{\D} (db) + \Sigma \D(r_i) b^{n-i} = 0$.

Due to separability $f'(b) \neq 0$. This implies $$\tilde{\D} (db) = - \Sigma \D(r_i) b^{n-i} / f'(b), 1 \leq i \leq n.$$

Since $c \in R[b], c = \Sigma s_i b^i, s_i \in R, 0 \leq i \leq n-1$. Due to the fact that $f(b) = 0$ and $\D(r_i) \in R$ for $1 \leq i \leq n$, we obtain
 $$\tilde{\D} (cdb) = c \tilde{\D} (db) = \Sigma \mu_i b^i/f'(b),  \mu_i  \in  R, 0 \leq i \leq n-1.$$
 Since $\Tr(b^i/f'(b)) = 0$ for $i < n-1$ and $\Tr (b^{n-1}/f'(b)) = 1$ ([13]), we have $\L(cdb)=\Sigma \mu_i \Tr (b^i/f'(b)) \in R$  \dotfill (12)

From (11) and (12) $\L (\omega)  \in  R$ for any $\omega  \in  \Omega_{B/k}$ and thus the claim is proved.

\smallskip
{\bf 3.2 Remarks.} Notations are as in the proof of the above theorem.

\smallskip
{\bf 1.} We have $\Omega_{B/k} = \underrightarrow{\lim} \Omega_{C/k}$ and $\Gamma_{B/k} = \underrightarrow{\lim} \Gamma_{C/k}$. Since $C$ is complete d.v.r. containing a field $k$ of characteristic 0, $C$ is geometrically regular over $k$. It follows by theorem 2.3 that $\Omega_{C/k}$ is a flat torsion free $C$-module and $\Gamma_{C/k} = 0$. Hence $\Omega_{B/k}$ is a flat torsion free $B$-module and $\Gamma_{B/k} = 0$.

\smallskip
{\bf 2.} We have $\Omega_{B/R} = \underrightarrow{\lim} \Omega_{C/R}$. Since $\Omega_{C/R}$ is a cyclic $C$-module of finite length, $\H{_m^0}(\Omega_{B/R}) = \Omega_{B/R}$ and by lemma 2.1 $\H{_m^1}(\Omega_{B/R}) = 0$.

\smallskip
{\bf 3.3 Corollary.} Notations are as in the proof of the above theorem.

\smallskip
The following sequence $0 \rightarrow \Omega_{B/R} \rightarrow \Omega_{R/k} \otimes \H{_m^1} (B) \rightarrow \H{_m^1} (\Omega_{B/k}) \rightarrow 0$ is exact.

From Remark 1 and the exact sequence: $0 \rightarrow \Omega_{R/k} \otimes B \rightarrow \Omega_{B/k} \rightarrow \Omega_{B/R} \rightarrow 0$ we obtain the following exact sequence:

$0 \rightarrow \H{_m^0} (\Omega_{B/R}) \rightarrow \H{_m^1}(\Omega_{R/k} \otimes B) \rightarrow \H{_m^1} (\Omega_{B/k}) \rightarrow 0$.

Since $\Omega_{R/k}$ is a flat $R$-module, we have $\H{_m^1} (\Omega_{R/k} \otimes B) \simeq \Omega_{R/k} \otimes \H{_m^1}(B)$.

The assertion now follows from Remark 2.

\smallskip
\textbf{3.4.} Next we deal with the non-complete d.v.r. case.

\textbf{Proposition.} Let $(R, m)$ be a discreet valuation ring containing a field $k$ of characteristic $0$ with a non-null derivation $\D \in \Der_k (R)$ and let $B$ be its integral closure in an algebraic extension $F$ of $Q(R)$. Then $\Hom_R(\Omega_{B/k}, R) \neq 0$.

\smallskip
{\bf Proof.}
Let $K$ denote the field of fractions of $R$. Since characteristic of $k=0$, $\D$ can be extended uniquly to a derivation $\D$ (same notation): $B \rightarrow F$. Let $\Tr: F \rightarrow K$ denote the trace map. We define an $R$-linear map $\L: \Omega_{B/k} \rightarrow K$ in the following way:

for any $\omega \ \in \ \Omega_{B/k}, \L(\omega) = \Tr(\tilde{\D}(\omega))$, where $\tilde{\D}:\Omega_{B/k} \rightarrow F$ is the $B$-linear map induced by $D$.

\smallskip
{\bf Claim.} Im $\L \subset R$.

Can assume $\tilde{\D} (\omega) \neq 0$. We will show that $\L (\omega) \in R$.
Let $C$ be a normal domain integral over $R$ such that $[Q(C): K] < \infty$ and $C \subset B$. We have $\Omega_{B/k} = \underrightarrow{\lim} \ \Omega_{C/k}$ over all such $C$ as above. For any $R \subset C \subset C'$, $C$, $C'$ as above are module finite extensions of $R$ and hence are semi-local Dedekind domains. Since characteristics of $k = 0$, for any such $C, k \hookrightarrow C$ is geometrically regular and hence $\Omega_{C/k}$ is a flat $C$-module (theorem 2.3). The exact sequences (5)-(8) in 3.1 are valid here. In particular we have the following exact sequence:

\smallskip
$0 \rightarrow \Omega_{C/k} \otimes C' \rightarrow \Omega_{C'/k} \rightarrow \Omega_{C'/C} \rightarrow 0$ \dotfill (13)

\smallskip

The left exactness in (13) follows from lemma 2.2 and the JZ sequence 2.5. Without any lose of generality can assume $\omega \in \Omega_{C/k}$ for some $C$ as above. Let $Q_1, \dots, Q_r$ be the maximal ideals of $C$ lying over $m$.
Then $\widehat{C} \simeq \widehat{C}_{Q_1} \times \dots \times \widehat{C}_{Q_r}$ (m-adic completion) \dotfill (14)

Here each $\widehat{C}_{Q_i}$ is a complete d.v.r. and a module finite extension of $\widehat{R}$. Since $R$ contains a field of characteristic 0, each $C_{Q_i} \rightarrow \widehat{C}_{Q_i}$ is geometrically regular and hence so is $C \rightarrow \widehat{C}$. By theorem 2.3 \ \ \ \ $\Omega_{\widehat{C}/C}$ is a flat $\widehat{C}$-module and $\Gamma_{\widehat{C}/C} = 0$. Hence, due to the JZ sequence 2.5, we obtain the following short exact sequence:

$0 \rightarrow \Omega_{C/k} \otimes \widehat{C} \rightarrow \Omega_{\widehat{C}/k} \rightarrow \Omega_{\widehat{C}/C} \rightarrow 0$ \dotfill (15)

We also have $\Omega_{\widehat{C}/k} \simeq \Omega_{{\widehat{C}_{Q_1}}/k} \times \dots \times \Omega_{{\widehat{C}_{Q_r}}/k}$ \dotfill (16)

Due to the fact that $\Omega_{B/k} = \underrightarrow{\lim} \ \Omega_{C/k}$ and (13), (15) and (16) we can write

$\omega = \omega_1 \times \dots \times \omega_r, \ \omega_i \in \Omega_{\widehat{C}_{{Q_i}/k}}, 1 \leq i \leq r$. \dotfill (17)

Since $\D \in \Der_k (R)$ and $\Der_k (R, \widehat{R}) \simeq \Der_k(\widehat{R})$ ($R$ being local any such $\D$ is continuous), $\D$ can be extended uniquely to a derivation $\widehat{\D}:\widehat{R} \rightarrow \widehat{R}$. This derivation $\widehat{\D}$ can be extended uniquely to a derivation $\widehat{\D}_i:\widehat{C}_{Q_i} \rightarrow K_i$, where $K_i$ is the field of fractions of $\widehat{C}_{Q_i}, 1 \leq i \leq r$. Thus $\D$ can be extended uniquely to a derivation $\D': \widehat{C} \rightarrow K_1 \times \dots \times K_r$ such that $\D'|\widehat{C}_{Q_i} = \widehat{\D}_i$. Due to uniqueness of extension of derivation we have $\D'|C = D$.

Hence from (17) we have $\tilde{\D}(\omega) = \Sigma\tilde{\widehat{\D}}_i (\omega_i)$ where $\tilde{\widehat{\D}}_i$ is the corresponding $\widehat{C}_{Q_i}$-linear map: $\Omega_{\widehat{C}_{Q_i}/k} \rightarrow K_i$. This implies that $L (\omega) = \Tr (\tilde{\D}(\omega)) = \Sigma \Tr (\tilde{\widehat{\D}}_i (\omega_i))$. By theorem 3.1, $\Tr (\tilde{\widehat{\D}}_i (\omega_i)) \ \in \ \widehat{R}, 1 \leq i \leq r$.

Hence $L(\omega) \ \in K \ \cap \widehat{R} = R$ and the claim is proved.

\smallskip
\textbf{Remark.} Remarks in 3.2 and corollary of 3.3 are also valid here.

\smallskip
{\bf 3.5.} Our next theorem in this section addresses the general set-up.

\smallskip
\textbf{Theorem.} Let $R$ be a normal domain containing a field $k$ of characteristic $0$ with a non-null derivation $\D \in \Der_k (R)$ and let $B$ be its integral closure in an algebraic extension $F$ of $Q(R)$. Then $\Hom_R(\Omega_{B/k}, R) \\ \neq 0$. In particular if $(R, m)$  is a complete local normal domain of dimension $n$ containing a field $k$ of characteristic 0 with a derivation as above, then $\H{_m^n} (\Omega_{B/k}) \neq 0$.

\smallskip
{\bf Proof.}
Let $K$ denote the field of fractions of $R$. Since characteristic of $k=0$, $\D$ can be extended to a derivation $\D$ (same notation): $B \rightarrow F$ and hence to a derivation: $B_P \rightarrow F$ for any prime ideal $P$ of $R$. Let $\Tr: F \rightarrow K$ denote the trace map. We define an $R$-linear map $\L: \Omega_{B/k} \rightarrow K$ in the following way:

for any $\omega \in \Omega_{B/k}, \L(\omega) = \Tr(\tilde{\D}(\omega))$ where $\tilde{\D}:\Omega_{B/k} \rightarrow F$ is the $B$-linear map induced by $\D$; the same prescription works for $\Omega_{{B_P}/k} \rightarrow F$ for any prime ideal $P$ of $R$.

\smallskip
{\bf Claim.} Im $\L \subset R$.

Can assume $\tilde{\D} (\omega) \neq 0$. Let $P$ be a prime ideal of height 1 in $R$. Let $\omega$ also denote the image of $\omega$ in $\Omega_{{B_P}/k}$. It follows from proposition 3.4 that  $\L(\omega) \in R_P$.

Since this holds for every prime ideal $P$ of height 1 in $R$ and $R$ is a normal domain, $\L(\omega) \in \underset{htP=1}{\cap}R_P = R$. Hence the claim follows.

In particular if $R$ is a complete local normal domain then the canonical module of $R$ is isomorphic to a height 1 ideal of $R$. Hence the assertion follows by local duality.

\smallskip
\textbf{3.6.} The following theorem deals with complete local domains. We recall that if $(C, q)$ is a complete local domain of dimension $n$ with coefficient ring $V(k)$, then $\Omega_{C/V} (\Omega_{C/k})$ is not necessarily a finitely generated $C$-module.

\smallskip
\textbf{Theorem.} Let $(C, q)$ be a complete local domain of dimension $n \geq 2$ with coefficient ring $V$ (field $k$). Let $(R, m)$ be a power series ring over $V(k)$ contained in $C$ such that $C$ is a module finite extension of $R$. We have the following: \\

i) if $Q(C)/Q(R)$ is separably algebraic, then rank $\widehat{\Omega}_{C/k} = n$, rank $\widehat{\Omega}_{C/V}=n-1$ and $\H_q^n(\Omega_{C/V}) (\H_q^n(\Omega_{C/k})) \neq 0$, \\

ii) if characteristic of $k =p >0, k$ perfect and $\Omega_{Q(C)/k} \neq 0$, then rank $\widehat{\Omega}_{C/k}=[Q(C):Q(C)^p]$, $\H_q^n(\Omega_{C/k}) \neq 0$ and  \\

iii) if $Q(C)/Q(R)$ is separably algebraic, then $\H_q^n(\Omega_{C/R}) = 0$; moreover if $\Omega_{C/R} \neq 0$ and $C$ is normal, then $\H_q^{n-1}(\Omega_{C/R}) \neq 0$.

\smallskip

\textbf{Proof.} Let $R = k[[X_1, \dots, X_n]]$ or $R = V[[X_1, \dots, X_{n-1}]]$ be contained in $C$ such that $C$ is a module finite extension of $R$. Let $m$ denote the maximal ideal of $R$. Then for any $C$-module $M, \H_m^i(M) = \H_q^i(M)$ for $i \geq 0$. We intend to show that $\widehat\Omega_{C/V} (\widehat\Omega_{C/k})$, the $m$-adic, equivalently $q$-adic, completion of $\Omega_{C/V} (\Omega_{C/k})$, is a finitely generated $C$-module of finite rank and hence a finitely generated $R$-module of finite rank. This would imply that $\Hom_R(\widehat\Omega_{C/V} (\widehat\Omega_{C/k}), R) \neq 0$. Since $R$ is complete we have
$$\Hom_R (\Omega_{C/V} (\Omega_{C/k}), R) = \Hom_R(\widehat\Omega_{C/V} (\widehat\Omega_{C/k}), R) \neq 0$$
Now the result would follow by local duality.

i) a) Equicharacteristic case. The arguments in (21.9.5, [5]) show that $\widehat\Omega_{C/k}$ is a finitely generated $C$-module and rank$_C  \ \widehat\Omega_{C/K} = n$.

\smallskip
b) Mixed characteristic p. In this case $C$ is a module finite extension of $R = V [[X_1, \dots, X_{n-1}]]$. We need to modify the arguments in (21.9.5, [5]). For this purpose we need the following lemma.

\textbf{Lemma.} Let $(R, m)$ be as above. Then $\Omega_{R/V}$ is a faithfully flat $R$-module, $\widehat{\Omega}_{R/V} \simeq R^{n-1}$ and the sequence
$$ 0 \rightarrow \underset{t\geq0}{\cap m^t} \Omega_{R/V} \rightarrow \Omega_{R/V} \rightarrow \widehat{\Omega}_{R/V} \rightarrow 0$$ is split exact.

Proof of the Lemma. There exists an exact sequence:

$0 \rightarrow R/m \rightarrow m/m^2 \rightarrow \Omega_{R/V} \otimes R/m \rightarrow 0$

\ \ \ \ \ \ \ \ \ \ \ $\overline{1} \rightarrow \overline{p}$

This implies that $\dm \ (\Omega_{R/V} \otimes R/m) = n-1$. Hence, due to theorem 2.3, $\Omega_{R/V}$ is a faithfully flat $R$-module. Moreover $T$, the Hausdorff-module associated with $\Omega_{R/V}$, i.e. $\Omega_{R/V}/ \cap m^t \Omega_{R/V}$, is generated by images of $dX_1, \dots, dX_{n-1}$ over $R$. $R$ being complete implies $T$ is also complete and hence $T = \widehat\Omega_{R/V}$. Due to the existence of derivations $\partial/\partial X_i : R \rightarrow R$, $\partial/\partial X_i (X_j) = \delta_{ij}, 1 \leq i, j \leq n-1$, and duality between derivatives and differentials it follows that images of $dX_1, \dots, dX_{n-1}$ are linearly independent in $\widehat\Omega_{R/V}$. Thus $\widehat\Omega_{R/V}$ is free $R$-module of rank $n-1$. And this implies the split exactness of the sequence in the lemma.

Now arguing as in (21.9.5, [5]) we observe that $\widehat\Omega_{C/V}$ is a finitely generated $C$-module of finite rank $n-1$ and hence a finitely generated $R$-module of finite rank.

\smallskip
ii) Characteristic of $k = p > 0$. We assume that $k$ is perfect and $\Omega_{Q(C)/k} \neq 0$. Let $R_1 = k[[X{_1^p}, \dots, X{_n^p}]]$. By (21.9.4, [5]) $\widehat \Omega_{C/k} = \Omega_{C/R_1}$. let $A = R_1[C^p]$. Then $\Omega_{C/R_1} = \Omega_{C/A}$. Let $y_1, \dots, y_t$ in $C$ form a $p$-basis of $Q(C)$ over $Q(R_1)$ i.e. $Q(C)$ over $k$, since $k$ is perfect. Due to our assumption that $\Omega_{Q(C)/k} \neq 0$ such a basis exists. Let $B = A[y_1, \dots, y_t]$. Then $B$ is a free A-module of rank $p^t$ and $\Omega_{B/A}$ is a free B-module of rank t with basis $dy_1, \dots, dy_t$. We have the following short exact sequence:
$$\Omega_{B/A} \otimes C \rightarrow \Omega_{C/A} \rightarrow \Omega_{C/B} \rightarrow 0$$
Since $Q(C) = Q(B)$, $\Omega_{C/B}$ is a torsion C-module.

\textbf{Claim.} $dy_1, \dots, dy_t$ are also linearly independent in $\Omega_{C/A}$.

The Claim follows via the existence of linearly independent $Q(A)$ derivations $\partial/\partial y_i, Q(C) \rightarrow Q(C), \partial/\partial y_i (y_j) = \delta_{ij}, 1 \leq i, j \leq t$.

Due to the above exact sequence we conclude that $\Omega_{C/A}$ is a finitely generated $C$-module of rank $t$. Hence $\widehat \Omega_{C/k}$ is a finitely generated $R$-module of finite rank.

The referee made the following observation: if $Q(C)/Q(R)$ is not separable, then since $R \rightarrow C$ s module-finite,
$$\Omega_{C/R} \otimes Q(C) \simeq \Omega_{C/R} \otimes Q(C) \simeq \Omega_{Q(C)/Q(R)} \neq 0.$$
Since $\widehat{\Omega}_{C/k} \rightarrow \widehat{\Omega}_{C/R}$ is onto, it follows that $\widehat{\Omega}_{C/k} \otimes Q(C) \neq 0$ and  thus $\widehat{\Omega}_{C/k}$ has a positive rank.

\smallskip
iii) Now assume that $R$ is a complete regular local ring contained in $C$ such that $C$ is a module finite extension of $R$ and $Q(C)/Q(R)$ is separably algebraic. Then by corollary to lemma 2.1, $H_m^n(\Omega_{C/R}) =0$.

Assume in addition that $\Omega_{C/R} \neq 0$ and $C$ is normal.

Let $0 \rightarrow R \rightarrow K \xrightarrow{\phi} \underset{htP=1}{\oplus} E(k(P)) \rightarrow \dots \rightarrow E \rightarrow 0$ denote an injective resolution of $R$, $K = Q(R)$. Let $T =$ Im$\phi$; then $T \hookrightarrow \underset{htP=1} {\oplus} E(k(P))$ is an essential extension.

Consider the exact sequence: $0 \rightarrow R \rightarrow K \rightarrow T \rightarrow 0$. By assumption Hom$_R(\Omega_{C/R}, K) =$ Hom$_K(\Omega_{(Q(C)/K}, K) = 0$ and Ext$_R^i (\Omega_{C/R}, K) =0$ for every $ i > 0$. Hence Hom$_R(\Omega_{C/R}, T) \simeq$ Ext$_R^1(\Omega_{C/R}, R)$. Thus, by local duality, in order to prove our assertion we need to show that Hom$_R(\Omega_{C/R}, T) \neq 0$.

We have an injection: $$\Hom_R(\Omega_{C/R}, T) \hookrightarrow \Hom_R(\Omega_{C/R}, \underset{htP=1}{\oplus} E(k(P))).$$

If Hom$_R(\Omega_{C/R}, E(k(P))) =$ Hom$_{R_P}(\Omega_{C_P/R_P}, E(k(P))) = 0$ for every prime ideal $P$ of $R$ of height 1, then $\Omega_{C_P/R_P} = 0$ for every prime ideal of height 1. Hence $\Omega_{C_q/R} = 0$ for every prime ideal $q$ of height 1 in $C$. By the theorem on purity of branch locus ([8], [13]) this would imply $\Omega_{C/R} = 0$, a contradiction. Thus there exists at least one prime ideal $P$ of height 1 such that Hom$_R(\Omega_{C/R}, E(k(P))) \neq 0$. Hence Hom$_R(\Omega_{C/R}, \underset{htP=1}{\oplus} E(k(P))) \neq 0$. \\ Let $f (\neq 0) \in$ Hom$_R(\Omega_{C/R}, \underset{htP=1}{\oplus} E(k(P)))$. $\Omega_{C/R}$ is a finitely generated $R$-module. Let $\omega_1, \dots, \omega_t$ denote the generators of $\Omega_{C/R}$; then each $f(\omega_i), 1 \leq i \leq t$ has only finitely many non-null components in $\underset {htP=1}{\oplus} E(k(P))$. Since $T \hookrightarrow \underset{htP=1} {\oplus} E(k(P))$ is essential, there exists $s \neq 0$ in $R$ such that $sf(\omega_i) \in T$ for $1 \leq i \leq t$ and $sf(\omega_i) \neq 0$ for at least one $i$. Thus Hom$_R(\Omega_{C/R}, T) \neq 0$.

\smallskip
{\bf Remarks. 1)} The proof for part i) of the above theorem is valid in equicharacteristic when dimension $C = 1$.

{\bf 2)} With hypothesis as in part iii) of the above theorem, dim$\Omega_{C/R} = n-1$. The proof in this part is also valid in the non-complete case.

\smallskip
\textbf{3.7.} Our final result of this section can be viewed as an extension of the theorem on purity of branch locus ([10], [15]) to infinitely generated integrally closed domains $B$ integral over a regular local ring $R$. Let us recall that an extension $R \rightarrow S$ is unramified if and only if $\Omega_{S/R} = 0$ and a local extension $(R, m) \rightarrow (S, q)$ of essentially finite type is unramified if and only if $mS =q$ and $S/q$ is a separable algebraic extension of $R/m$.

\smallskip
\textbf{Theorem.} Let $R$ be a regular local ring and let $B$ denote the integral closure of $R$ in a field extension $K$ of $Q(R)$ such that $K/Q(R)$ is separably algebraic. Suppose that for every prime ideal $q$ of height 1 in $B$, $\Omega_{{B_q}/R} = 0$. Then $\Omega_{B/R} = 0$ i.e. $B$ is an unramified extension of $R$.

\smallskip
\textbf{Proof.} If $[K: Q(R)] < \infty$ then this is the usual theorem on purity of branch locus. Assume that $[K: Q(R)] = \infty$. Let $R \subset C \subset C' \subset B$ be such that $C, C'$ are integral closures of $R$ in $Q(C), Q(C')$ respectively and $[Q(C): Q(R)] < \infty$, $[Q(C') : Q(R)] < \infty$. Then $C, C'$ are module-finite extensions of $R$ and $C'$ is a module finite extension of $C$. Since $B = \underrightarrow{\lim} C, \ \Omega_{B/R} = \underrightarrow{\lim} \Omega_{C/R}$ over all such $C$ as mentioned above. Let $q$ be a prime ideal of $B$ of height 1 and let $P = q \cap R$. Then $\Omega_{{B_P}/R} = \underrightarrow{\lim} \Omega_{{C_P}/R}$. From the assumption it follows that $\Omega_{{B_P}/R} = 0$. Each $C_P$ is a semi-local Dedekind domain and $C_P \hookrightarrow C'_P$ is a module finite extension. Hence $\Gamma_{{C'_P}/{C_P}} = 0$ (lemma 2.2). By the JZ sequence 2.5 we obtain an exact sequence

$0 \rightarrow \Omega_{{C_P}/R} \otimes C'_P \rightarrow \Omega_{{C'_P}/R} \rightarrow \Omega_{{C'_P}/{C_P}} \rightarrow 0$

Since $\Omega_{{B_P}/R} =0$, from the above exact sequence it follows that $\Omega_{{C_P}/R} = 0$ for every such $C$. By our assumption this holds for every prime ideal $P$ of height 1 in $R$ and hence for prime ideals $\tau$ of height 1 in $C$ we have $\Omega_{C_{\tau/R}} = 0$. Due to the theorem on purity of branch locus it follows that $\Omega_{C/R} = 0$ for every such $C$. Hence $\Omega_{B/R} = 0$.

\smallskip
\section*{Section 4: Mixed characteristic} \label{s2}

In this section we deal with the situation $R \hookrightarrow B$ where $B$ is the absolute integral closure of $R$. In the positive characteristic $p > 0$ case we have $\Omega_{B/R} = 0$, $\Omega_{B/k} = 0$, i.e. $k \rightarrow B$, $R \rightarrow B$ are unramified extensions. Due to theorem 3.5 we now concentrate on the mixed characteristic case. First we define formally unramified extension (19.10, [5]).

\smallskip
\textbf{Definition.} Let $f:A \rightarrow B$ be a continuous ring homomorphism between two topological rings. We say $B$ is formally unramified over $A$ if given any discrete $A$-algebra $C$, any ideal $I$ of $C$ such that $I^2 = 0, C/I$ with discrete topology and $\eta:C \rightarrow C/I$, the natural surjection, every continuous $A$-algebra homomorphism $\phi:B \rightarrow C/I$ has at most one $A$-algebra extension $\psi:B \rightarrow C$ such that $\phi = \eta . \psi$.

\smallskip
\textbf{4.1. Theorem.} Let $(R, m)$ be a complete regular local ring of dimension n in mixed characteristic $p$ with coefficient ring $V$ and let $B$ denote its absolute integral closure. We have the following:

{\bf i)} $\Omega_{B/V} = p \ \Omega_{B/V}, \Omega_{B/R} = p \ \Omega_{B/R}$. These observations imply the following:

\smallskip
a) For any ideal $I$ of $R$ of the form $I=(p, x_1, .., x_{i-1})$ generated by a part of a system of parameters containing $p$ of $R$, $H_I^i(\Omega_{B/V})=0$ and $H_I^i(\Omega_{B/R})=0$. In particular $H_m^n(\Omega_{B/V})=0$ and $H_m^n(\Omega_{B/R})=0$. 

\smallskip
b) $\widehat{\Omega}_{B/V}=0=\widehat{\Omega}_{B/R}$, where $(\widehat{-})$ denotes the I-adic completion. And hence $V \rightarrow B, R \rightarrow B, V \rightarrow \widehat{B}$ and $R \rightarrow \widehat{B}$ are all formally unramified extensions.

\smallskip
{\bf ii)} Non-vanishing of $\H{_m^{n-1}}(\Omega_{B/R})$ implies the direct summand property for integral extensions of $R$ ([5]) (equivalently descent of flatness for integral extensions of Noetherian rings ([10])) and the converse is also true.

Non-vanishing of $\H{_m^{n-1}} (\Omega_{B/V})$ implies the direct summand property for integral extensions of $R$.

\smallskip
\textbf{Proof.}
{\bf i)} Actually we would prove the following more general statement:

Let $A$ be an integral domain of characteristic 0 and let $B$ be the absolute integral closure of an integral domain $C$ containing $A$. Let $t$ be a positive integer such that $t$ is not a unit in $C$. Then for any proper ideal $I$ of $A$ containing $t$ and for every $n > 0$ we have $\Omega_{B/A} = I^n \Omega_{B/A}$. Let $\widehat{\Omega}_{B/A}$ and $\widehat{B}$ denote the $I$-adic completion of $\Omega_{B/A}$ and $B$ respectively. Then $\widehat{\Omega}_{B/A} = 0 = \widehat{\Omega}_{\widehat{B}/A}$ and $A \rightarrow B, A \rightarrow \widehat{B}$ are formally unramified extensions with respect to the $I$-adic topology.

\smallskip
\text{Proof.} For every element $x(\neq 0)$ in $B$, $x$ has a $t^nth$ root $y$ in $B$ i.e. $y^{t^n} = x$. Hence $dx = t^ny^{t^n-1}dy$; this implies that $\Omega_{B/A} = t^n \Omega_{B/A}$ and thereby $\Omega_{B/A} = I^n \Omega_{B/A}$ for every $n >0$.

Let $\tilde B = B/I^nB$, $\tilde A = A/I^nA$; via base change we obtain
$$ \Omega_{B/A} \otimes \tilde B \simeq \Omega_{\tilde B/\tilde A}$$
We have from above, $\Omega_{\tilde B / \tilde A} = 0$.
Hence $\widehat \Omega_{B/A} = 0$. Let $\tilde{\widehat{B}} = \widehat{B}/{I^n \widehat{B}} = B/I^nB$; this implies $\Omega_{\tilde{\widehat{B}}/\tilde A} = \Omega_{\tilde B/\tilde A} = 0$. Hence $\widehat \Omega_{\widehat B/A} = 0$. It follows from (20.7.4 and 20.7.5, [5]) that $A \rightarrow B, A \rightarrow \widehat{B}$ are formally unramified extensions.

\smallskip
Next let $(A, q)$ be a complete local normal domain of dimension $n$ with coefficient ring $V$ and let $B$ be its absolute integral closure. By our assumption it follows from above that $\Omega_{B/V} = p^t \Omega_{B/V}$ and if $I = (p, x_1, .., x_{i-1})$ generated by a  part of a system of parameters of $R$ and $I_t = (p^t, x_1^t, .., x_{i-1}^t)$, then $\Omega_{B/V} = I_t \Omega_{B/V}$ for $t>0$.  Similar statement is valid for $\Omega_{B/R}$. Now the assertion follows immediately.

\smallskip
{\bf ii)} Let $\alpha: R \hookrightarrow B$ denote the inclusion as mentioned in the setup. We need to show that $\alpha$ splits as an $R$-module map $\Leftrightarrow \H{_m^{n-1}} (\Omega_{B/R}) \neq 0$.

Recall that $\alpha$ splits as an $R$-module map $\Leftrightarrow B^* = \Hom_R (B, R) \neq 0$ i.e. $\H{_m^{n}} (B) \neq 0$ ([8]) (same condition for descent of flatness ([12])).

First we deal with the unramified case i.e. $R = V[[X_1, \dots, X_{n-1}]]$. If dimension of $R =1$, $R = V, \Omega_{Q(B)/Q(R)}=0$ implies that $\Omega_{B/R} [1/p] = 0$; hence $\H{_m^0} (\Omega_{B/R}) = \Omega_{B/R} \neq 0$.

We can assume $\dm R \geq 2$.

Since the direct summand property is valid over $R$ ([2]), we have $\H{_m^n} (B) \neq 0$. Since $V \rightarrow R$ is geometrically regular, $\Omega_{R/V}$ is a flat $R$-module (theorem 2.3).

By lemma 2.4 we have the following short exact sequence:
$$ 0 \rightarrow \Omega_{R/V} \otimes_R B \rightarrow \Omega_{B/V} \rightarrow \Omega_{B/R} \rightarrow 0$$

Since $\Omega_{R/V}$ is a flat $R$-module, $\H{_m^i} (\Omega_{R/V} \otimes B) \simeq \Omega_{R/V} \otimes \H{_m^i} (B)$ for every $i > 0$.

Applying local cohomology to the above short exact sequence, by i), we obtain the following short exact sequence:

\smallskip
$\rightarrow \H{_m^{n-1}} (\Omega_{B/R}) \rightarrow \Omega_{R/V} \otimes_R \H{_m^n} (B) \rightarrow 0$. \dotfill (18)

\smallskip
\ \ \ \ \ \ \ \ \ \ \ \ Due to the short exact sequence

$0 \rightarrow R/m \rightarrow m/m^2 \rightarrow \Omega_{R/V} \otimes R/m \rightarrow 0$

\ \ \ \ \ \ \ \ \ \ $(\overline{1} \rightarrow \overline{p})$

we obtain $\Omega_{R/V} \otimes R/m \neq 0$ and thus $\Omega_{R/V}$ is faithfully flat over $R$.
Since $\H{_m^{n}} (B) \neq 0$, $\Omega_{R/V} \otimes_R \H{_m^{n}}(B) \neq 0$. Hence from (18) it follows that $\H{_m^{n-1}} (\Omega_{B/R}) \neq 0$.

\smallskip
Conversely, suppose that $\H{_m^{n-1}}(\Omega_{B/R}) \neq 0$. Consider the short exact sequence obtained via base change: $R \hookrightarrow B$ (spectral sequence)
$$0\rightarrow \Ext{_B^{1}}(\Omega_{B/R}, B^*)\rightarrow \Ext{^1_{R}}(\Omega_{B/R}, R)\rightarrow \Hom_B(\Omega_{B/R}, \Ext{_R^{1}}(B, R)).$$
Note that since $(R, m)$ is complete, given any $R$-module $M$, Ext$_R^i (M, R)$ is $m$-separated - actually it is linearly compact for every $i \geq 0$ ([9]). This can be checked by considering a free resolution $F_\bullet$ of $M$ and then applying Hom$_R (F_\bullet, R)$ and noticing that for any continuous $R$-linear map $f: T_1 \rightarrow T_2$ between two linearly compact $m$-separated $R$-modules Im$f$ and Ker$f$ are both linearly compact, $m$-separated closed submodules of $T_2$ and $T_1$ respectively ([9]).

Since $\Ext{_R^{1}}(B, R)$ is $m$-separated and $\Omega _{B/R} = p^n \Omega_{B/R}$ for every $n > 0$, the right end term of the above exact sequence is 0. By local duality, $\H{_m^{n-1}}(\Omega_{B/R}) \neq 0$ implies $\Ext{_R^{1}} (\Omega_{B/R}, R) \neq 0$. Hence from the above exact sequence $\Ext{_B^{1}}(\Omega _{B/R}, B^*) \neq 0$ and thus $B^* \neq 0$.

If $\H{_m^{n-1}} (\Omega_{B/V}) \neq 0$, same argument as above shows that $\Hom_R (B, R) \\ \neq 0$.

Now let us deal with the ramified case. In this case the proof that $\Ext{_R^1}(\Omega_{B/R}, R) \neq 0$ implies $\Hom_R (B, R) \neq 0$ is same as above. Our proof for the reverse implication needs a modification due to the fact that in this case $V \rightarrow R$ is not geometrically regular and hence $\Omega_{R/V}$ is not necessarily flat. We have the following exact sequence:

$ 0\rightarrow N \rightarrow \Omega_{R/V} \otimes_R B \xrightarrow{\phi} \Omega_{B/V} \rightarrow \Omega_{B/R} \rightarrow 0$  \dotfill (19)

Since $R$ is ramified, $R = S[X]/(f(X))$, where $S = V[[X_1, \dots, X_{n-1}]]$ and $f(X)$ is an Eisenstein polynomial in $S[X]$. Then $\Omega_{R/S} = R/(f')$. By lemma 2.2 and the JZ sequence we have the following short exact sequence:

$ 0 \rightarrow \Omega_{S/V} \otimes_S R \rightarrow \Omega_{R/V} \rightarrow {R/f'R} \rightarrow 0$

Tensoring this sequence with $B$ we obtain the following exact sequence:

$ 0 \rightarrow \Omega_{S/V} \otimes_S B \rightarrow \Omega_{R/V} \otimes_R B \rightarrow {B/f'B} \rightarrow 0$  \dotfill (20)

\ \ \ \ \ \ \ \ \ Let $q$ denote the maximal ideal of $S$. Then for any $R$-module $M, H_q^i (M) = H_m^i (M)$ and $\Hom_R (B, R) \simeq \Hom_S (B, S)$. Since $\H{_m^n} (B) \neq 0$ and $\Omega_{S/V}$ is faithfully flat over $S$, it follows that $\H{_m^n} (\Omega_{S/V} \otimes B) \simeq \Omega_{S/V} \otimes \H_m^n(B) \neq 0$. We have $(\Omega_{S/V} \otimes \H_m^n (B))^\vee \simeq \Hom(\Omega_{S/V}, \H_m^n(B)^\vee) \simeq \Hom (\Omega_{S/V}, B^*) \simeq \underset{1}{\overset{n-1}{\oplus}} B^*$ (lemma in  i) b) theorem 3.6). Hence from the above exact sequence we obtain $\H{_m^n} (\Omega_{R/V} \otimes B) \neq 0$.

We also obtain an exact sequence (lemma 2.4) via $V \rightarrow S \rightarrow B$:

$ 0 \rightarrow \Omega_{S/V} \otimes_S B \rightarrow \Omega_{B/V} \rightarrow \Omega_{B/S} \rightarrow 0$ \dotfill (21)

In (19) let $W = \im \phi$. From (19), (20) and (21) we get the following commutative diagram of exact sequences:

\[
\begin{CD}
@. @. 0 @. 0 \\
@. @. @VVV  @VVV \\
@. @. N @= N @. \\
@. @. @VVV  @VVV \\
0 @>>>\Omega_{S/V} \otimes B @>>> \Omega_{R/V} \otimes B @>>> B/f'B @>>> 0 \\
@.   \parallel @. @VVV @VVV @. \dots \dots (22)\\
0 @>>>\Omega_{S/V} \otimes B @>>> W @>>> U @>>>0 \\
@. @. @VVV  @VVV \\
@. @. 0 @. 0
\end{CD}
\]

From the exact column on the right end we obtain $\H{_m^n} (N) = 0$. Hence from the middle column we obtain $\H{_m^n} (\Omega_{R/V} \otimes B) \simeq \H{_m^n} (W) \neq 0$. Now applying local cohomology to the exact sequence (from 19):

$$ 0 \rightarrow W \rightarrow \Omega_{B/V} \rightarrow \Omega_{B/R} \rightarrow 0$$

We obtain our required result $\H{_m^{n-1}} (\Omega_{B/R}) \neq 0$ (by part i).

When $\Omega_{S/V}$ is not flat over $S$, i.e. when $S=V$, $R=S[X]/(f(X))$ as mentioned earlier, $\dim R = 1$ and $\dim B=1$. Then $\Omega_{B/V}$ is p-torsion and $H_p^0 (\Omega_{B/V}) = \Omega_{B/V} \neq 0$.

\smallskip
\textbf{Corollaries.} With $R, B$ as in the above theorem we have the following

\smallskip
i) $\H{_m^{n-1}} (\Omega_{B/R}) \neq 0$.

This follows due to Andr\'e's proof ([2]) of the direct summand conjecture.

ii) $\Ext^1_B (\Omega_{B/R}, B^*) \neq 0$.

This follows from the proof of the theorem.

\smallskip
\textbf{Remark.} We do not know whether part ii) of the above theorem is valid in equicharacteristic zero.

\smallskip
\textbf{4.2.} Our next proposition points out a relation between local cohomologies of $B$ and that of $\Omega_{B/R}, \Omega_{B/V}$.

\smallskip
\textbf{Proposition.} Let $(R, m)$ be a complete unramified regular local ring of dimension $n$ with coefficient ring $V$ and let $B$ be its absolute integral closure.

\smallskip
Then, for every $i \geq 0$, there exists a short exact sequence:
$$0\rightarrow \H{_m^i}(\Omega_{B/V}) \rightarrow \H{_m^i}(\Omega_{B/R})\rightarrow \Omega_{R/V} \otimes_R \H{_m^{i+1}} (B)\rightarrow 0$$

\textbf{Proof.} Consider the exact sequence (lemma 2.4):
$$0\rightarrow \Omega_{R/V} \otimes_RB\rightarrow \Omega_{B/V} \rightarrow \Omega_{B/R} \rightarrow 0$$. Since $\Omega_{R/V}$ is $R$-flat and $\Omega_{R/V} \otimes_R{R/m} \neq 0, \Omega_{R/V}$ is faithfully flat as an $R$-module. Hence $\H{_m^i}(\Omega_{R/V} \otimes_RB) = \Omega_{R/V} \otimes_R \H{_m^i}(B)$ is non-null whenever $\H{_m^i}(B)$ is so.

\medskip
\textbf{Claim.} For any $B$-module $M$, $\Omega_{B/V} \otimes_B \H{_m^i}(M) = 0$. Same result holds for $\Omega_{B/R}$.

\medskip
\textbf{Proof of the Claim.} $\Hom_R (\Omega_{B/V} \otimes_B \H{_m^i} (M), E) \simeq \Hom_B (\Omega_{B/V}, \\ \Hom_R(\H{_m^i}(M), E)) \simeq \Hom_B(\Omega_{B/V}, \Ext{_R^{n-i}} (M, R))$ (local duality). Since $\Ext{_R^{n-i}} (M, R)$ is Hausdorff in the $m$-adic topology and $\Omega_{B/V} = p^n \Omega_{B/V}$, for $n>0$, we have $\Hom_B(\Omega_{B/V}, \Ext{_R^{n-i}} (M, R)) = 0$ and hence the claim follows.

\smallskip
Now the assertion follows from noticing that the map: $\H{_m^i} (\Omega_{R/V} \otimes_R B) \rightarrow \H{_m^i} (\Omega_{B/V})$ is the composition of the following two maps:

\smallskip
$\Omega_{R/V} \otimes_R \H{_m^i} (B) \simeq \Omega_{R/V} \otimes_R B \otimes_B \H^i_m(B) \rightarrow \Omega_{B/V} \otimes_B \H{_m^i}(B)$ and $\Omega_{B/V} \otimes_B \H{_m^i}(B)\rightarrow \H{_m^i}(\Omega_{B/V})$

\smallskip
When $\Omega_{R/V}$ is not $R$-flat, i.e. when $\dim R=1$ i.e. $R=V, \Omega_{R/V} =0$ and the asserted sequence becomes obvious.

\medskip
\textbf{Corollary 1.} With $R, B$ as above, $\H{_m^{i+1}} (B) \neq 0 \Rightarrow \H{_m^i} (\Omega_{B/R}) \neq 0$ for $i \geq 0$.

\smallskip
\textbf{Corollary 2.} With $R, B$ as above, for every $i \geq 0$ there exists a short exact sequence:

$0 \rightarrow \underset{1}{\overset{n-1}{\oplus}} \Ext {^{n-i-1}_R} (B, R) \rightarrow \Ext{_R^{n-i}} (\Omega_{B/R}, R) \rightarrow \Ext{_R^{n-i}} (\Omega_{B/V}, R) \rightarrow 0$.

In particular

$0 \rightarrow \underset{1}{\overset{n-1}{\oplus}} B^* \rightarrow \Ext{_R^1} (\Omega_{B/R}, R) \rightarrow \Ext{_R^1} (\Omega_{B/V}, R) \rightarrow 0$
is exact.

\smallskip
\textbf{Proof.} It follows from the lemma in part b), i) theorem 3.6 that $\Omega_{R/V}$ is a faithfully flat $R$-module, $\hat{\Omega}_{R/V} \simeq R^{n-1}$ and the sequence
$$ 0 \rightarrow \underset{t\geq0}{\cap m^t} \Omega_{R/V} \rightarrow \Omega_{R/V} \rightarrow \hat{\Omega}_{R/V} \rightarrow 0$$ is split exact. Since $R$ is complete, $\Ext{_R^{n-i-1}} (B, R)$ is m-separated (a proof of this fact is sketched in the proof of part ii), theorem 4.1). We have $(\Omega_{R/V} \otimes \H{_m^{i+1}} (B))^\vee \simeq \Hom_R (\Omega_{R/V}, \H{_m^{i+1}} (B)^\vee ) \simeq \Hom_R (\Omega_{R/V}, \\ \Ext{^{n-i-1}_R} (B, R))$. Since Ext$_R^{n-i-1} (B, R)$ is m-separated, it follows from the above lemma that $\Hom_R (\Omega_{R/V}, \Ext{^{n-i-1}_R} (B, R)) \simeq \Hom_R (\underset{1} {\overset{n-1}{\oplus}} R, \\ \Ext{^{n-i-1}_R} (B, R)) \simeq \underset{1} {\overset{n-1}{\oplus}} \Ext{^{n-i-1}_R} (B, R)$. Now applying $\Hom_R (-, E)$ to our asserted exact sequence in the above proposition the corollary follows.

\bigskip

\textbf{Remarks.}

\smallskip
{\bf 1.} If $\dm R \geq 2$, then $\H{_m^0} (\Omega_{B/V}) \simeq \H{ _m^0} (\Omega_{B/R})$. 

This follows due to the fact that if $\dim R \geq 2$, then $H_m^i (B)=0, i=0, 1$ and $\Omega_{R/V}$ is $R$-flat.

\smallskip
{\bf 2.} Part i) of theorem 4.1 (except $\H_m^n (\Omega_{B/R}) = 0$) and results in 4.2 are not necessarily valid in equicharacteristic zero.

This follows mainly due to theorem 3.5.

\bigskip
\textbf{4.3.} Our final theorem deals with local cohomology of module of differentials with respect to the ideal ($x_1, \dots, x_{n-1}$) generated by a part of a regular system of parameters.

 \textbf{Theorem.} Let $(R, m)$ be a complete regular local ring of dimension $n$ in mixed characteristics $p$ with co-efficient ring $V$ and let $B$ be its absolute integral closure. Let $p, x_1, \dots, x_{n-1}$ be a system of parameters of $R$ such that $x_1, \dots, x_{n-1}$ form a regular system of parameters. Let $\underline{x}$ denote the ideal $(x_1, \dots, x_{n-1})$ in $R$. Then $\H{_{\underline{x}}^{n-1}}(\Omega_{B/V}) \neq 0$. Non-vanishing of $\H{_{\underline{x}}^{n-1}}(\Omega_{B/R})$ implies the direct summand property for integral extensions of $R$.
\smallskip

\textbf{Proof.} For any $x \neq 0 \in R$ and for any $R$-module $M$, let $\tilde{K} (x; M)$ denote the complex: $0 \rightarrow M \rightarrow M [1/x] \rightarrow 0$  and $\tilde{K} (x_1, \dots, x_{n-1}; R)$ denote the complex $\underset{1}{\overset{n-1}{\otimes}} \tilde{K} (x_i; R)$. The local cohomology modules of $\Omega_{B/V}$ with respect to $m$ are obtained by taking cohomologies of \\ $\tilde{K} (x_1, \dots, x_{n-1}; R) \otimes \tilde{K} (p; \Omega_{B/V})$. Since this complex is a tensor product of two complexes, we obtain the following exact sequence:

\smallskip
$\H{_m^{n-1}} (\Omega_{B/V}) \rightarrow \H{_{\underline{x}}^{n-1}} (\Omega_{B/V}) \xrightarrow{\eta} \H{_{\underline{x}}^{n-1}} (\Omega_{B/V} [ 1/p]) \rightarrow \H{_m^n} (\Omega_{B/V}) \rightarrow 0 .........................................................................(23)$

\smallskip
By theorem 4.1 $\H{_m^n} (\Omega_{B/V}) = 0$; hence $\eta$ is onto.

Let $q$ denote the prime ideal $\underline{x} \subset R$, then $q R_q$ is a maximal ideal in $R_q$. We have:

$\H{_{\underline{x}}^{n-1}} (\Omega_{B/V} [ 1/p]) = \Omega_{B/V}[1/p] \underset{R}{\otimes} \H{_{\underline{x}}^{n-1}} (R) \simeq \Omega_{B/V} \underset{R}{\otimes} \H{_{\underline{x}}^{n-1}} (R[1/p]) \simeq \Omega_{B/V} \underset{R}{\otimes} \H{_{qR_q}^{n-1}} (R_q) \simeq \Omega_{B/V} \underset{R}{\otimes} \H{_{q\widehat{R}_q}^{n-1}} (\widehat{R}_q) \simeq \H{_{q\widehat{R}_q}^{n-1}} (\Omega_{{B_q}/k} \otimes \widehat{R}_q)$, ($k$ = field of fractions of $V$). By theorem 3.5 we have $\Hom_{R_q}(\Omega_{B_{q/k}}, R_q) \neq 0 \Rightarrow \Hom_{R_q}(\Omega_{B_{q/k}}, \widehat{R}_q) \neq 0 \Rightarrow \Hom_{\widehat{R}_q}(\Omega_{B_{q/k}} \otimes \widehat{R}_q, \widehat{R}_q) \neq 0$. Hence, by local duality, it follows from above that $\H{_{\underline{x}}^{n-1}} (\Omega_{B/V} [ 1/p]) \neq 0$.

Since $\eta$ is onto in (23) it follows that $\H_{\underline{x}}^{n-1} (\Omega_{B/V}) \neq 0$.

Now let us assume $\H_{\underline{x}}^{n-1} (\Omega_{B/R}) \neq 0$. By similar arguments as for exactness in (23) we obtain the following exact sequence

\smallskip
$\H_m^{n-1} (\Omega_{B/R}) \rightarrow \H_{\underline{x}}^{n-1} (\Omega_{B/R}) \rightarrow \H_{\underline{x}}^{n-1} (\Omega_{B/R} [ 1/p]) \rightarrow \H_m^n (\Omega_{B/R}) \rightarrow 0 ...................................................................(24)$

\smallskip
By arguments as above we have $\H_{\underline{x}}^{n-1} (\Omega_{B/R} [1/p]) = \H_{qRq}^{n-1}(\Omega_{B_q/R_q})=0$ (corollary to lemma 2.1) and hence $\H_{\underline{x}}^{n-1} (\Omega_{B/R}[1/p]) = 0$. It follows  from (24) that the sequence $\H_m^{n-1} (\Omega_{B/R}) \rightarrow \H_{\underline{x}}^{n-1} (\Omega_{B/R}) \rightarrow 0$ is exact. Hence the assertion follows by part ii) theorem 4.1.

\smallskip
\textbf{Corollary.} If $\H{_p^0} (\H{_{\underline{x}}^{n-1}} (\Omega_{B/V})) \neq 0$, then $\H{_m^{n-1}} (\Omega_{B/V}) \neq 0$.

Suppose $\H{_p^0} (\H{_{\underline{x}}^{n-1}} (\Omega_{B/V})) \neq 0$. Then, in the exact sequence (23), $\eta$ is not an isomorphism and hence $\H{_m^{n-1}} (\Omega_{B/V}) \neq 0$.

\smallskip
\textbf{Remark.} With notations as above if $I$ is an ideal in $R$ generated by a system of parameters of length $n-1$ such that $p \in I$, then $\H_I^{n-1} (\Omega_{B/R}) = 0$ and $\H_I^{n-1} (\Omega_{B/V}) = 0$. 

This has been pointed out in part i) of theorem 4.1.

\smallskip
\textbf{4.4.} We here mention a couple of questions that came up during this study and whose answers we do not know at present.

1. Let $(R, m)$ be a complete regular local ring of dimension $n$ of equicharacteristic $0$ and let $B$ be its absolute integral closure. Is $\H_m^{n-1}(\Omega_{B/R}) \neq 0$?

\medskip
2. Let $(R, m)$ be a complete regular local ring of dimension $n$ of mixed characteristic $p$ with coefficient ring $V$ and let $B$ be its absolute integral closure. Let $p, x_1, \dots, x_{n-1}$ be a system of parameters of $R$ such that $x_1, \dots, x_{n-1}$ form a part of a regular system of parameters.

\smallskip
a) Is $\H_{\underline{x}}^{n-1}(\Omega_{B/R}) \neq 0$?

\smallskip
b) Is $\H_p^0 (\H_{\underline{x}}^{n-1}(\Omega_{B/V})) \neq 0$?

\smallskip
\textbf{Remark.} It can be shown that an affirmative answer for b) implies an affirmative answer for a).

\end{document}